\documentclass[12pt]{amsart}

\usepackage[utf8]{inputenc}
\usepackage[margin=0.8in]{geometry}
\usepackage{amsmath}
\usepackage{amsthm}
\usepackage{amssymb}
\usepackage{amsrefs}
\usepackage{mathrsfs}
\usepackage[normalem]{ulem}
\usepackage{enumitem}
\usepackage[bottom,flushmargin]{footmisc}
\usepackage[mmddyyyy]{datetime}
\usepackage{gensymb}
\usepackage{mathtools}
\usepackage{graphicx}
\usepackage{setspace}
\usepackage{wasysym}
\usepackage{bbm}
\usepackage{ulem}

\allowdisplaybreaks

\usepackage{color}   
\usepackage{hyperref}
\hypersetup{
	colorlinks=true, 
	linktoc=all,     
	linkcolor=blue,  
}

\usepackage{tikz}
\usetikzlibrary{shapes,positioning,intersections,quotes}
\usepackage[font=small,labelfont=bf]{caption}

\DeclareFontFamily{U}{mathx}{\hyphenchar\font45}
\DeclareFontShape{U}{mathx}{m}{n}{
	<5> <6> <7> <8> <9> <10>
	<10.95> <12> <14.4> <17.28> <20.74> <24.88>
	mathx10
}{}
\DeclareSymbolFont{mathx}{U}{mathx}{m}{n}
\DeclareFontSubstitution{U}{mathx}{m}{n}
\DeclareMathAccent{\widecheck}{0}{mathx}{"71}
\DeclareMathAccent{\wideparen}{0}{mathx}{"75}

\makeatletter
\newcommand*\rel@kern[1]{\kern#1\dimexpr\macc@kerna}
\newcommand*\widebar[1]{%
	\begingroup
	\def\mathaccent##1##2{%
		\rel@kern{0.8}%
		\overline{\rel@kern{-0.8}\macc@nucleus\rel@kern{0.2}}%
		\rel@kern{-0.2}%
	}%
	\macc@depth\@ne
	\let\math@bgroup\@empty \let\math@egroup\macc@set@skewchar
	\mathsurround\z@ \frozen@everymath{\mathgroup\macc@group\relax}%
	\macc@set@skewchar\relax
	\let\mathaccentV\macc@nested@a
	\macc@nested@a\relax111{#1}%
	\endgroup
}
\makeatother

\newcommand{\prs}[1]{\left(#1\right)}
\newcommand{\prss}[1]{(#1)}
\newcommand{\sbrks}[1]{\left[#1\right]}
\newcommand{\abs}[1]{\left|#1\right|}

\newcommand{\cbrks}[1]{\left\{#1\right\}}
\newcommand{\lprs}[1]{\left(#1\right.}

\newcommand{\lsbrks}[1]{\left[#1\right.}

\newcommand{\rprs}[1]{\left.#1\right)}

\newcommand{\rsbrks}[1]{\left.#1\right]}

\DeclarePairedDelimiter{\floor}{\lfloor}{\rfloor}

\newcommand{\norm}[1]{\left\lVert#1\right\rVert}
\newcommand{\norms}[1]{\lVert#1\rVert}

\newcommand{\lskip}{\vspace{0.4cm}}
\newcommand{\slskip}{\vspace{0.2cm}}

\newcommand{\nlspace}{\\[0.3cm]}

\newcommand\numberthis{\addtocounter{equation}{1}\tag{\theequation}}

\newcommand{\rnm}{\mathbb{R}}
\newcommand{\znm}{\mathbb{Z}}
\newcommand{\cnm}{\mathbb{C}}

\newcommand{\eun}[1]{\rnm^{#1}}

\newcommand{\ra}{\rightarrow}

\newcommand{\se}{\subseteq}

\newcommand{\finv}{^{-1}}

\DeclareMathOperator{\sspan}{span}

\DeclareMathOperator{\spt}{\text{sppt}}

\newcommand{\fconst}{2\pi i}

\newcommand{\mbm}[1]{\mathbbm{#1}}
\newcommand{\mc}[1]{\mathcal{#1}}

\newcommand{\os}[1]{^{#1}}

\newcommand{\wh}[1]{\widehat{#1}}

\newcommand{\nf}{\infty}
\newcommand{\tm}{\times}
\newcommand{\pd}{\partial}
\newcommand{\dx}[1]{\,d#1}

\newcommand{\cl}{\colon}

\newcommand{\tx}[1]{\text{#1}}

\newtheorem{theorem}{Theorem}
\newtheorem{corollary}[theorem]{Corollary}

\newtheorem{lemma}[theorem]{Lemma}
\newtheorem{proposition}[theorem]{Proposition}
\newtheorem{heuristic}[theorem]{Rule of Thumb}

\theoremstyle{definition}

\newtheorem{remark}{Remark}

\newcommand{\ft}{\mc{F}}

\renewcommand{\phi}{\varphi}

\newcommand{\sph}[1]{S\os{#1}}
\newcommand{\cosp}[1]{\cos\prs{#1}}

\newcommand{\chr}[1]{\mbm{1}_{#1}}

\DeclareMathOperator{\conv}{\tx{Conv}}

\newcommand{\subskip}{0.2in}

\begin{document}
	
\title[The Triangle Operator]{The Triangle Averaging Operator}
\author{Eyvindur A. Palsson}
\address{Department of Mathematics, Virginia Tech, Blacksburg, VA 24061, USA}
\email{palsson@vt.edu}
\thanks{The work of the first listed author was supported in part by Simons Foundation Grant \#360560. The authors would like to thank Danqing He for alerting us to a gap in the argument of a previous version of this work.}
\author{Sean R. Sovine}
\address{Department of Mathematics, Virginia Tech, Blacksburg, VA 24061, USA}
\email{sovine5@vt.edu}

\maketitle

\begin{abstract}
	We examine the averaging operator $T$ corresponding to the manifold in $\eun{2d}$ of pairs of points $(u,v)$ satisfying $\abs{u} = \abs{v} = \abs{u - v} = 1$, so that $\{0,u,v\}$ is the set of vertices of an equilateral triangle. We establish $L^p \times L^q \ra L^r$ boundedness for $T$ for $(1/p, 1/q, 1/r)$ in the convex hull of the set of points $\{\prs{0, 0, 0},\, \prs{ 1, 0 , 1 }, \prs{ 0, 1, 1 }, \prss{ {1}/{p_d}, {1}/{p_d}, {2}/{p_d} }\}$, where $p_d = \frac{19d - 4}{11d - 12}$ and $d \geq 7$. 
\end{abstract}


\section{Introduction}

\vskip.125in

Much work in harmonic analysis has centered around determining the mapping properties and regularity of operators given by the convolution $f \star \mu$ of a function $f$ and a measure $\mu$ and maximal variants of these operators. Often this work has addressed the case where $\mu$ is supported on a submanifold satisfying a curvature condition, as in \cite{GreenleafPrincipal}, \cite{SteinWainger}, \cite{SoggeSteinI}. Operators of this type have applications in continuous geometric combinatorics, including applications to generalizations of the Falconer distance problem; see \cite{FalconerDistance}, \cite{IosevichIncidence}, \cite{GreIos}, \cite{GrafMulGen}, for example. There has also been some interest in multilinear versions of surface averages, i.e., multilinear convolution operators of the form
\begin{align}\label{eqn8}
	(f_1,\ldots, f_n) \mapsto [(f_1\otimes\cdots\otimes f_n)\star \mu] (x, \ldots, x),
\end{align}
where $f_1, \ldots, f_n$ are measurable functions on $\eun{d}$ and $\mu$ is a Borel measure on $\eun{nd}$.
As an early example of such operators, in \cite{OberlinMultiConv} Oberlin introduced the multilinear convolution of $d$ functions on $\rnm$ with the unit sphere in $\eun{d}$ and completely characterized the boundedness of this operator. 
In \cite{GrafakosSoria} Grafakos and Soria proved some fundamental and general results regarding bilinear convolutions with positive measures. 
Grafakos, Greenleaf, Iosevich, and Palsson \cite{GrafMulGen} proved general results on applications of bounds on multilinear generalized Radon transforms to point configuration problems. A multilinear generalized radon transform is essentially a version of \eqref{eqn8} in which the submanifold that $\mu$ is supported on is allowed to vary depending on~$x$.

In Geba, et al.\ \cite{GebaEtAl}, the bilinear version of the spherical maximal operator, given by
\begin{align*}
\mc{B}(f,g)(x) = \sup_{t > 0}\abs{\int_{\sph{2d-1}}f(x - tu)g(x - tv)\dx{\sigma(u,v)}},
\end{align*}
was studied as a basic example in a class of  bilinear maximal averaging operators corresponding to averages with respect to finite Borel measures with well-behaved Fourier transforms. 
The study of this operator was later taken up by Barrionuevo, Grafakos, He, Honzik, and Oliveira in the paper \cite{Barrio}. These authors used a square function technique similar to that used by Rubio de Francia \cite{Rubio} to bound the linear maximal spherical averaging operator, along with a wavelet decomposition of the Fourier transform of the spherical surface measure, to establish boundedness on $L^p \tm L^p \ra L^{p/2}$ for a particular exponent $p$ with $p/2 < 1$.
This bound was then interpolated against bounds for the range of exponents with $1 < p,q,r \leq \nf$ obtained by taking the $L^\nf$ norm of one input function and applying results on linear Fourier multiplier operators.
Most recently, Jeong and Lee \cite{JeongLee} established a complete characterization of the $L^p \tm L^q \ra L^r$ boundedness of the bilinear spherical maximal operator $\mc{B}$ by using the slicing identity 
\begin{align*}
\int_{\sph{2d-1}} F(x,y)\dx{\sigma(x,y)} = \int_{B^d(0,1)}\int_{\sph{d-1}}F\prs{x, \sqrt{1 - \abs{x}^2}}(1 - \abs{x}^2)\os{\frac{d-2}{2}}\dx{\sigma_{d-1}(y)}\dx{x}
\end{align*}
to majorize $\mc{B}(f,g)(x)$ pointwise by a product of the linear spherical maximal operator applied to one input function and the Hardy-Littlewood maximal operator applied to the other input function. These authors also obtain Lorentz space estimates for endpoint cases. 

Greenleaf and Iosevich \cite{GreIos} introduced the bilinear convolution operator
\begin{equation*}
	B(f,g)(x) = \int f(x - u)g(x - v)\dx{K(u,v)}
\end{equation*}
where $dK$ is the surface measure on the manifold $\{ (u,v)\in \eun{2}\tm \eun{2} ~\cl~ \abs{u} = \abs{v} = \abs{u - v} = 1 \}$ of pairs of points $(u,v)$ such that $\{0, u, v\}$ is the set of vertices of an equilateral triangle in the plane with side length 1. In \cite{GreIos} the authors proved the estimates
\begin{align*}
	\norm{B(f,g)}_{L^1(\eun{2})} \leq C\norm{f}_{L^2_{-\beta_1}(\eun{2})}\norm{g}_{L^2_{-\beta_2}(\eun{2})} ~~ \tx{ if }~~ \beta_1 + \beta_2 = \frac{1}{2}, ~~ \beta_1, \beta_2 \geq 0,
\end{align*}
for positive functions $f$ and $g$. Using these estimates the authors proved that, when $E\se \eun{2}$ is a compact set with Hausdorff dimension greater than $\frac{7}{4}$, the set of three-point configurations determined by $E$ has positive Lebesgue measure as a subset of $\eun{3}$, where a three-point configuration $C = (x,y,z)\in E\tm E \tm E$ is identified with its triple of pairwise distances. This is a variation of the Falconer distance problem, in which it is conjectured that for a compact set $F\se \eun{d}$ of Hausdorff dimension $\dim_{\mc{H}}(F) > \frac{d}{2}$ the one-dimensional Lebesgue measure of the distance set $\Delta(F) = \{ \abs{x - y} ~\cl~ x, y \in F \}$ is positive. The Falconer distance problem is in turn a continuous version of the Erd\H{o}s distinct distance problem.

Here we study a generalization of the bilinear operator $B$ introduced by Greeneleaf and Iosevich in \cite{GreIos} that applies to functions on $\eun{d}$, for $d \geq 3$. We define 
\begin{equation*}
	T(f,g)(x) = \int_M f(x - u)g(x - v)\dx{\mu(u,v)}
\end{equation*}
to be the averaging operator corresponding to the surface measure $\mu$ on the submanifold of $\eun{2d}$
\begin{equation*}
	M = \{ (u,v) \in \eun{d}\tm \eun{d} ~\cl~  \abs{u} = \abs{v} = \abs{u - v} = 1\},
\end{equation*}
consisting of all pairs of points $(u,v)$ such that $\{0,u,v\}$ is the set of vertices of an equilateral triangle of side length 1 with one vertex at the origin.

We use interpolation and techniques from \cite{Barrio} to show that this operator is bounded from $L^p(\eun{d})\tm L^q(\eun{d}) \ra L^r(\eun{d})$ for a range of indices $p,q,r$ with $1\leq p,q \leq \nf$ that properly includes the Banach range, by which we mean the range for which $1/p + 1/q \leq 1$. 
Our proof requires that $d \geq 7$; the restriction that $d \geq 5$ cannot be avoided using the tools that we use in our proof. Our main result is:


\begin{theorem}
	For $d\geq 7$ the operator $T$ is bounded from
	\begin{equation*}
		L^p(\eun{d})\tm L^q(\eun{d}) \ra L^r(\eun{d})
	\end{equation*}
	for $\prs{\frac{1}{p}, \frac{1}{q}, \frac{1}{r}}$ in the convex hull of the set of points
	\begin{align*}
		\cbrks{  
			\prs{0, 0, 0},\, \prs{ 1, 0 , 1 }, \prs{ 0, 1, 1 }, \prs{ \frac{1}{p_d}, \frac{1}{p_d}, \frac{2}{p_d} }
		},
	\end{align*}
	where $p_d = \frac{19d - 4}{11d - 12}$. 
\end{theorem}

\begin{center}
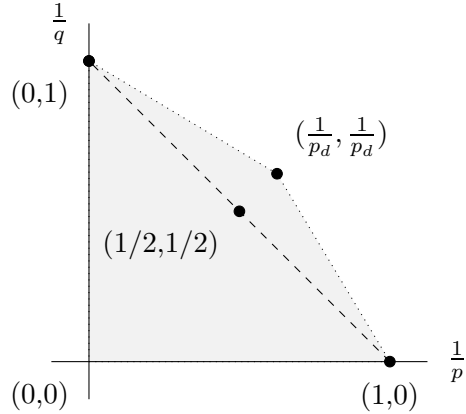

	\begin{tikzpicture}[scale=0.5]
	\draw[fill=gray!10, style = dotted] (0,0) -- (8,0) -- (40/8,40/8) -- (0,8) -- cycle;
	\draw[color=black] (-1,0) -- (9,0);
	\draw[color=black] (0,-1) -- (0,9);
	\draw[color=black, style=dashed] (8,0) -- (0,8);
	
	\node [fill, draw, circle, minimum width=4pt, inner sep=0pt] at (4,4) {};
	\node[below left=2pt of {(4,4)}, outer sep=2pt] {\small(1/2,1/2)};
	
	\node [fill, draw, circle, minimum width=4pt, inner sep=0pt] at (40/8,40/8) {};
	\node[above right=1pt of {(40/8,40/8)}, outer sep=2pt] {\small$(\frac{1}{p_d},\frac{1}{p_d})$};
	
	\node [fill, draw, circle, minimum width=4pt, inner sep=0pt] at (8,0) {};
	\node[below=2pt of {(8,0)}, outer sep=2pt] {\small(1,0)};
	
	\node [fill, draw, circle, minimum width=4pt, inner sep=0pt] at (0,8) {};
	\node[below left=2pt of {(0,8)}, outer sep=2pt] {\small(0,1)};
	
	\node [fill, draw, circle, minimum width=4pt, inner sep=0pt] at (0,8) {};
	\node[below left=2pt of {(0,0)}, outer sep=2pt] {\small(0,0)};
	
	\node[right=2pt of {(9,0)}, outer sep=2pt] {$\frac{1}{p}$};
	
	\node[left=2pt of {(0,9)}, outer sep=2pt] {$\frac{1}{q}$};
	\end{tikzpicture}
	\captionof{figure}{\color{black} Region of boundedness we obtain for $T$. }\label{fig1}
\end{center}
\vspace{.2in}


\begin{remark}
By analogy with the linear spherical averaging operator, we conjecture that $T$ is bounded on $L^p\tm L^q$ for $p, q \in [1,\nf]$. Our Theorem 1 does not obtain boundedness for this full range of exponents. This situation seems similar to that of  the bilinear spherical maximal operator, where Barrionuevo, et al.\ \cite{Barrio} were able to obtain a non-trivial range of boundedness using multiplier techniques, and Jeong and Lee \cite{JeongLee} obtained the full range of boundedness using a geometrically-motivated majorization. However, though they don't give the full range of boundedness, the multiplier techniques in \cite{Barrio} and the decomposition used in this work can be applied to other multipliers that satisfy similar decay estimates. One can verify that the decay estimate \eqref{eqn1} that we obtain for the mutiplier $\wh{\mu}$ is not sufficient to allow the application of the bilinear H\"ormander multiplier theorem (see Grafakos and Van Nguyen \cite{GrafakosVanNguyen}) or the the Bilinear Calder\'on-Vaillancourt theorem of Miyachi and Tomia \cite{MiyachiTomita} to this operator. 
\end{remark}

We also introduce and prove some initial results for the maximal version of $T$,
\begin{equation*}
\mc T(f,g)(x) := \sup_{t > 0}\abs{\int_M f(x - tu)g(x - tv)\dx{\mu(u,v)}}.
\end{equation*}
We show using a variant of the example given in \cite{Barrio} that if $\mc T$ is bounded from $L^p \tm L^q$ into $L^r$ with $1/p + 1/q = 1/r$, then $p, q > \frac{d}{d-1}$. By interpolation against bounds obtained by taking one of $p, q$ equal to $\infty$ we show that $\mc T$ is bounded on the half of this region in $(1/p,1/q)$-space for which $1/p + 1/q < \frac{d-1}{d}$. This is in contrast with the situation of the bilinear spherical maximal operator, which is bounded on the entire Banach range $1 < p, q, r \leq \nf$.

\subsection*{Novel contributions of this work:}
The novel contributions of this work are as follows: We describe the surface measure on $M$ in terms of integration with respect to the Haar measure on $SO(d)$, and this allows us to obtain an explicit representation of the Fourier transform $\wh{\mu}$ of the surface measure on $M$. Using this representation, we obtain an estimate of the decay of $\wh{\mu}(\xi,\eta)$, which is shown in \eqref{eqn1}. The decay in our estimate has significant non-uniformity, and in fact depends on the ratio $\abs{\xi}/\abs{\eta}$ and on the angle between $\xi$ and $\eta$ as vectors in $\eun{d}$. Finally, we apply a decomposition of $\wh{\mu}$ that is tailored to the decay properties given by \eqref{eqn1} and show that this decomposition allows us to establish the $L^p\tm L^q \ra L^r$ boundedness of $T$ for some exponents $p,q,r$ outside the Banach range. We also introduce the maximal version of the operator $T$ and establish by interpolation a range of $p,q$ for which this maximal operator is bounded on $L^p(\eun{d})\tm L^q(\eun{d})$, and we give an example that establishes lower bounds on the indices $p,q$ for which such boundedness can hold.

In summary, this work provides an interesting example of a Fourier multiplier with non-uniform decay. 
We make an initial attempt to address the boundedness of the corresponding multiplier operator using a suitable decomposition and some tools from the existing literature. We also make some initial progress towards characterizing the boundedness of the corresponding maximal operator. 


\subsection*{Overview of results:}
We compute the Fourier transform $\wh{\mu}$ below and obtain the estimate 
\begin{align}\label{eqn1}
\abs{\pd^\alpha_\xi \pd^\beta_\eta\wh{\mu}(\xi,\eta)} 
&\leq C_{\alpha, \beta}( 1 + \min\{\abs{\xi}, \abs{\eta}\}\abs{\sin\theta} )\os{-\frac{d-2}{2}}
\prs{ 1 + \abs{(\xi, \eta)}}\os{-\frac{d-2}{2}}.
\end{align}
The procedure we use to bound $T$ outside the Banach range is based on ideas from the proof in~\cite{Barrio}.
The main idea is to decompose the multiplier $m = \wh{\mu}$ dyadically away from the origin as $m = \sum_{i=0}^{\nf}m_i$, where $m_i$ is supported on a dyadic annular region on scale $\approx 2^i$, and then to obtain nicely decaying bounds for the operators corresponding to the $m_i$'s.
We obtain bounds $$\norm{T_i(f,g)}_{L^1}\leq C_i\norm{f}_{L^2}\norm{g}_{L^2}$$ for the operator $T_i$ corresponding to the multiplier $m_i$, with norms $C_i$ exponentially decreasing in~$i$, and we obtain bounds $$\norm{T_i(f,g)}_{L^{p/2}}\leq D_i \norm{f}_{L^p}\norm{g}_{L^p},$$ with norms $D_i$ increasing in $i$, using an estimate from \cite{Rubio}. We interpolate between these two bounds for $T_i$ as far as possible so that the interpolated norms $C_i\os{\theta}D_i\os{1 - \theta}$ still form a summable sequence in $i$, and summing these bounds for $T_i$ we obtain a bound for $T$, 
$$\norm{T(f,g)}_{L^{p/2}}\leq C \norm{f}_{L^p}\norm{g}_{L^p}$$
for a range of $p$ with $p/2 < 1$. 
When $p=\nf$ or $q = \nf$ we can majorize $T(f,g)(x)$ by the product of the $L^\nf$ norm of one input function and the spherical maximal function corresponding to the other function, and then we get a range of bounds for $T$ from the known boundedness of the linear spherical maximal operator. We then interpolate between these bounds and the trivial $L^\nf\tm L^\nf \ra L^\nf$ bound for $T$.

Most of the work in this proof lies in obtaining the bounds 
$$\norm{T_i(f,g)}_{L^1}\leq C_i\norm{f}_{L^2}\norm{g}_{L^2}$$
with exponentially decreasing norms $C_i$. We make essential use of the following proposition, which is a consequence of results proved by Grafakos, He, and Slav\'ikov\'a~\cite{GrafakosL2} using a wavelet technique.
\begin{proposition}\label{prop1}
	If $m \in L^q(\eun{2d}) \cap C_0^{M_1}(\eun{2d})$, where $M_1 = \floor*{\frac{2d}{3}} + 1$, and $m$ satisfies
	\begin{equation*}
	\norm{\pd^\alpha m}_{L^\nf} \leq C_0 < \nf,
	\end{equation*}
	for all $\abs{\alpha}\leq M_1$,
	then there is a constant $D$ depending only on $d$ and $q$ such that the bilinear operator $T_m$ with multiplier $m$ satisfies 
	\begin{align*}
	\norm{T_m}_{L^2\tm L^2 \ra L^1} \leq DC_0\abs{\mathrm{sppt}\, m}\os{\frac{1}{4}}.
	\end{align*}
\end{proposition}
\noindent Here $\abs{\mathrm{sppt}\, m}$ is the Lebesgue measure of the support of the function $m$. We cannot apply this lemma directly to the dyadic pieces $m_i$, because for these pieces we have only the uniform decay estimates
\begin{align*}
	\norm{\pd^\alpha m_i}_{L^\nf} \leq C2\os{-i\frac{d-2}{2}},
\end{align*}
and $\abs{\mathrm{sppt}\, m}\os{\frac{1}{4}} \approx \frac{id}{2}$. Thus direct application of Proposition \ref{prop1} gives an $L^2\tm L^2 \ra L^1$ norm for $T_i$ of $\approx 2\os{i}$. In order to obtain bounds for $T_i$ decaying in $i$ we need to further decompose $m_i$ to take advantage of the additional non-uniform decay in the estimate \eqref{eqn1} and then apply Proposition \ref{prop1} to the pieces of this finer decomposition. We decompose $m_i(\xi,\eta)$ dyadically based on the ratio $\abs{\xi}/\abs{\eta}$, and then further decompose these pieces dyadically based on the angle $\theta$ between $\xi$ and~$\eta$.

It is worth considering whether it is possible to obtain a geometric majorization of $T$, similar to the majorization of $\mc{B}$ utilized by Jeong and Lee \cite{JeongLee}. It follows from the description of $M$ below that we can express $T$ as 
\begin{align*}
	T(f,g)(x) &= \int_{\sph{d-1}}f(x - u)S_u(g)(x)\dx{\sigma_{d-1}(u)},
\end{align*}
where $S_u(g)(x)$ is the average of the function $g(x - \cdot)$ over the submanifold
\begin{equation*}
	N_u = \cbrks{ v \in \sph{d-1} ~\cl~ \abs{u - v} = 1 } = \sph{d-1}\cap \cbrks{ v ~\cl~ u\cdot v = 1/2 },
\end{equation*}
which is a $(d-2)$-dimensional sphere of radius $\frac{\sqrt{3}}{2}$ centered at $x - \frac{1}{2}u$, lying inside the hyperplane containing $x - \frac{1}{2}u$ with normal direction $u$. 
A natural approach to majorize $T$ is to
take the supremum over $u$ in the operator $S_u$, giving 
\begin{align*}
	\abs{T(f,g)(x)} \leq \abs{S(f)(x)}\cdot\sup_{u\in \sph{d-1}}\abs{S_u(g)(x)},
\end{align*}
where $S(f)(x)$ is the linear spherical maximal function for $f$. 
This leads to the study of the directional maximal operator 
\begin{align*}
	S_*(g)(x) := \sup_{u\in \sph{d-1}}\abs{S_u(g)(x)}.
\end{align*}
We plan to consider the boundedness of $S_*$, but do not address the boundedness of this operator in this work.

\vskip.25in


\section{Describing the Manifold}

\vskip.125in

Given $(u,v) \in \sph{d-1}\tm\sph{d-1}$, we compute
\begin{align*}
	\norm{u - v}^2 &= 2 - 2\cos\theta,
\end{align*}
where $\theta$ is the angle between $u$ and $v$. Hence we see that $\norm{u  -v} = 1$ exactly when $\theta = \pm \pi/3$. Thus, since rotations preserve angles, we can represent any $(u,v)\in M$ as 
\begin{align*}
	(u,v) = \prs{R{e_1}, R\sbrks{\frac{1}{2}e_1 + \frac{\sqrt{3}}{2}e_2}},
\end{align*}
where $e_i$ is the $i$th standard basis vector in $\eun{d}$ and $R$ is an appropriately chosen rotation.

%
%
%

This description of $M$ suggests that we can effectively integrate over the $M$ using the following formula
\begin{align*}
\int_M f(x,y)\dx{\mu(x,y)} = 
	\int_{SO(d)}f\prs{ Re_1,  R\sbrks{\frac{1}{2}e_1 + \frac{\sqrt{3}}{2}e_2}}\dx{R} = \int_{SO(d)}f\prs{R\sbrks{\frac{1}{2}e_1 + \frac{\sqrt{3}}{2}e_2}, Re_1}\dx{R},
\end{align*}
where $SO(d)$ is the special orthogonal group in dimension $d$ and $dR$ is the Haar measure on $SO(d)$, and we have used the invariance of the Haar measure. 
We will give an argument to justify this identification below.

Because $SO(d)$ is compact, $dR$ is both left- and right-invariant. For any closed subgroup $\mc{H}$ of $SO(d)$ there is an invariant Radon measure $d[S]$ on the quotient $SO(d)/\mc H$ such that the following quotient integral formula holds
\begin{align*}
	\int_{SO(d)}f(R)\dx{R} = \int_{SO(d)/ \mc H}\int_{\mc H}f(SR')\dx{R'}\dx{[S]},
\end{align*}
where $dR'$ is the Haar measure on the subgroup $\mc H$ and $S$ is any representative of $S\mc{H}$.
Since $SO(d)$ is compact its Haar measure is finite, and we assume that all measures are normalized to be probability measures.\footnote{Since the measures involved are all probability measures and the map $S\mapsto  \int_{\mc H}f(SR')\dx{R'}$ is constant on cosets in $SO(d)/\mc{H}$ we have the identity
$$\int_{SO(d)/ \mc H}\int_{\mc H}f(SR')\dx{R'}\dx{[S]} = \int_{SO(d)/ \mc H}\int_{\mc H}\int_{\mc H}f(STR')\dx{R'}\dx{T}\dx{[S]} = \int_{SO(d)}\int_{\mc H}f(RR')\dx{R'}\dx{R}.$$
} 
For a reference on Haar measure on locally compact groups, see Deitmar and Echterhoff \cite{Deitmar}, and for a reference on matrix Lie groups see Hall \cite{HallLie}. 

If $\mc H$ is the closed subgroup of rotations fixing $e_1$, then each coset $R\mc{H}$ is the set of all rotations taking $e_1$ to $Re_1$, and $\mc{H}$ is isomorphic to $SO(d-1)$. Using the quotient integral formula with this choice of $\mc{H}$ lets us express the integral of $f$ with respect to the surface $\mu$ measure on $M$
\begin{align*}
	\int_M f(x,y)\dx{\mu(x,y)} &= \int_{SO(d)}\int_{SO(d-1)}f\prs{ Re_1, R\sbrks{ \frac{1}{2}e_1 +  \frac{\sqrt{3}}{2}(0, R'e_1\os{d-1})} }\dx{R'}\dx{R}\nlspace
	&= \int_{SO(d)}\int_{SO(d-1)}f\prs{ R\sbrks{ \frac{1}{2}e_1 + \frac{\sqrt{3}}{2}(0,R'e_1\os{d-1}) }, Re_1 }\dx{R'}\dx{R},
\end{align*}
where $e_i^{d-1}$ is the $i$th standard basis vector in $\eun{d-1}$. 

Now consider the measure of a ball $B(x, r)$ in $M$ as a metric space with the metric inherited from the ambient space $\eun{2d}$. From our definition of the measure $\mu$ on $M$ we can see that $\mu$ is a uniformly distributed measure. That is, for any two points $x, y$ in $M$ and any $r > 0$, $\mu(B(x,r)) = \mu(B(y,r))$. This follows from the fact that we can map the first ball onto the second by a map of the form $(x,y)\mapsto (Ux, Uy)$, where $U \in SO(d)$, and such maps preserve distances.
Our measure is invariant under such transformations by construction.
By a result of Christensen~\cite{Christensen}, such a uniformly distributed measure is unique up to a constant multiple. One can verify using the inverse function theorem that,
since the equations defining $M$ are invariant under the linear isomorphisms $(x,y) \mapsto (Ux, Uy)$, with $U \in SO(d)$, the natural surface measure on $M$ is also invariant under the action of such isomorphisms, and hence is uniformly distributed. Thus our expression for the integral over $M$ recovers the normalized surface measure on $M$.
\\[-2mm]

\noindent\textbf{The triangle averaging operator.}
We define the \textit{bilinear triangle averaging operator} as
\begin{align*}
	T(f,g)(x) &:= \int_M f(x - u)g(x - v)\dx{\mu(u,v)}\nlspace
	&= \int_{SO(d)}\int_{SO(d-1)}f\prs{x- Re_1}g\prs{x- R\sbrks{ \frac{1}{2}e_1 +  \frac{\sqrt{3}}{2}(0, R'e_1\os{d-1})} }\dx{R'}\dx{R}\nlspace
	&= \int_{SO(d)}\int_{SO(d-1)}f\prs{x- R\sbrks{ \frac{1}{2}e_1 +  \frac{\sqrt{3}}{2}(0, R'e_1\os{d-1})} }g\prs{x- Re_1}\dx{R'}\dx{R}.
\end{align*}
for Schwartz functions $f$ and $g$ on $\eun{d}$. 


\vskip.25in

\section{Computing the Fourier Transform of the Surface Measure}

\vskip.125in

Here we address the case where $d\geq 3$.
We compute
\begin{align}\label{eqn5}
\wh{\mu}(\xi, \eta) &= \int_{SO(d)}\int_{SO(d-1)}\exp\sbrks{-\fconst \prs{\xi\cdot Re_1 + \eta \cdot R\sbrks{ \frac{1}{2}e_1 + \frac{\sqrt{3}}{2}\prs{0, R'e_1\os{d-1}} } } }\dx{R'}\dx{R}.
\end{align}
Now we can use the invariance of the Haar measure on $SO(d)$ to precompose $R$ with a rotation whose transpose takes $\xi$ to a multiple of $e_1$ and $\eta$ to a vector in $\sspan\{e_1, e_2\}$. 
This map can be chosen to take $\eta$ to $\abs{\eta}\cdot e_1$ and $\xi$ to $\abs{\xi}\cos\theta\cdot e_1 + \abs{\xi}\sin\theta\cdot e_2$, where $\theta$ is the angle between $\xi$ and~$\eta$. 
This lets us write 
\begin{align*}
\wh{\mu}(\xi, \eta) &= \int_{SO(d)}\int_{SO(d-1)}\exp\lprs{-\fconst \lsbrks{
		\prs{ \abs{\xi}\cos\theta + \frac{1}{2}\abs{\eta} }
		e_1\cdot Re_1 }}\nlspace
		&\qquad\qquad\qquad\qquad\qquad\qquad\qquad
		\rprs{+\abs{\xi}\sin\theta\, e_2\cdot Re_1+\rsbrks{\abs{\eta}e_1 \cdot {\frac{\sqrt{3}}{2}R\prs{0, R'e_1\os{d-1}} } } }\dx{R'}\dx{R}.
\end{align*}
Now we can pull out some factors and use orthogonality and the fact that $\wh{\sigma}_{d-2}$ is radial\footnote{Below we will abuse notation and write $\wh{\sigma}_{d-2}(x) = \wh{\sigma}_{d-2}(\abs{x})$.} to obtain for the inner integral\footnote{Here we are splitting the integral over $SO(d-1)$ over cosets of the subgroup preserving $e_1^{d-1}$, similar to above, since the integrand here is only affected by the action of $R'$ on $e_1^{d-1}$.}
\begin{align*}
&\int_{SO(d-1)}\exp\prs{ \abs{\eta}e_1 \cdot R{\frac{\sqrt{3}}{2}\prs{0, R'e_1\os{d-1}} }}\dx{R'}\nlspace
%
%
=~& \wh{\sigma}_{d-2}\prs{ \frac{\sqrt{3}}{2}\abs{\eta}\cdot P_{-1}R^T e_1 },
\end{align*}
where $P_{-1}(u_1,\ldots, u_d)\mapsto (u_2,\ldots, u_d)$ is the projection on to the last $d-1$ coordinates, and $\sigma_{d-2}$ is the surface measure on the unit sphere in $\eun{d-1}$. However, if $P_1$ is the projection onto the first coordinate, then by the Pythagorean theorem we have
\begin{align*}
\abs{ P_{-1}R^T e_1 } = \sqrt{1 - \prs{ P_1R^Te_1 }^2} = \sqrt{ 1 - (R e_1\cdot e_1)^2 },
\end{align*}
so the inner integral is 
\begin{align*}
= \wh{\sigma}_{d-2}\prs{ \frac{\sqrt{3}}{2}\abs{\eta}\sqrt{ 1 - (R e_1\cdot e_1)^2  }}.
\end{align*}
Since the rotation $R'$ has been eliminated, we can rewrite the Fourier transform in terms of a spherical integral,
\begin{alignat*}{1}
\mkern-18mu\wh{\mu}(\xi, \eta) &= \int_{\sph{d-1}}
\exp\prs{-\fconst \sbrks{
		\prs{ \abs{\xi}\cos\theta + \frac{1}{2}\abs{\eta} }
		e_1\cdot v + \abs{\xi}\sin\theta\, e_2\cdot v}}
\wh{\sigma}_{d-2}\prs{ \frac{\sqrt{3}}{2}\abs{\eta}\sqrt{ 1 - (v\cdot e_1)^2  }}\dx{\sigma(v)}.\\[-4mm]
\end{alignat*}

Next we take note of the following simple slicing formula for the spherical integral,
\begin{align*}
\int_{\sph{d-1}}f(u)du &= \sum_{\pm}\int_{B_{d-1}(0,1)} f\prs{\pm\sqrt{1 - \abs{y}^2}, y}\frac{dy}{\sqrt{1 - \abs{y}^2}}\nlspace
&= \sum_{\pm}\int_{0}^{1}\int_{\sph{d-2}}f\prs{\pm\sqrt{1 - r^2}, r\omega}\frac{r\os{d-2}}{\sqrt{1 - r^2}}\dx{\sigma_{d-2}(\omega)}\dx{r}.\numberthis\label{eqn3}
%
\end{align*}
We can use this formula to take advantage of the fact that the integrand in the last expression above for $\wh{\mu}(\xi, \eta)$ depends only on the first two components of $y$. Applying this formula to $\wh{\mu}$ gives
\begin{align*}
&\wh{\mu}(\xi, \eta) =\nlspace
&\:\:\:\sum_{\pm}\int_{0}^{1}\int_{\sph{d-2}}
\exp\prs{-\fconst \sbrks{
		\pm\sqrt{1 - r^2}\prs{ \abs{\xi}\cos\theta + \frac{1}{2}\abs{\eta} }
		+ \abs{\xi}\sin\theta\, e_1^{d-1}\cdot r\omega}}\nlspace
&\qquad\qquad\qquad\qquad\qquad\qquad\qquad\qquad\qquad\qquad\qquad
\cdot
\wh{\sigma}_{d-2}\prs{ \frac{\sqrt{3}}{2}\abs{\eta}r}
\frac{r\os{d-2}}{\sqrt{1 - r^2}}\dx{\sigma_{d-2}(\omega)}\dx{r}\nlspace
%
%
&\:=2\int_{0}^{1}\cosp{ 2\pi \sqrt{1 - r^2}\prs{ \abs{\xi}\cos\theta + \frac{1}{2}\abs{\eta} } }
\wh{\sigma}_{d-2}\prs{ r\abs{\xi}\abs{\sin\theta} }
\wh{\sigma}_{d-2}\prs{ \frac{\sqrt{3}}{2}\abs{\eta}r}
\frac{r\os{d-2}}{\sqrt{1 - r^2}}\dx{r}.
\end{align*}
Inserting the well-known formula for $\wh{\sigma}_{d-2}$, this is 
\begin{align}\label{eqn2}
&\wh{\mu}(\xi,\eta) =
2(2\pi)^2\int_{0}^{1}\cosp{ 2\pi \sqrt{1 - r^2}\prs{ \abs{\xi}\cos\theta + \frac{1}{2}\abs{\eta} } }
\frac{J_{\frac{d-3}{2}}\prs{ 2\pi r \abs{\xi}\abs{\sin\theta} }}
{\prs{ r\abs{\xi}\abs{\sin\theta} }\os{\frac{d-3}{2}}}
\frac{ J_{\frac{d-3}{2}}\prss{2\pi r\frac{\sqrt{3}}{2}\abs{\eta} } }{\prss{ r\frac{\sqrt{3}}{2}\abs{\eta} }\os{\frac{d-3}{2}}}
\frac{r\os{d-2}}{\sqrt{1 - r^2}}\dx{r}.
\end{align}


\vskip.25in

\subsection*{Estimating Derivatives of $\wh{\mu}$}

To use Corollary \ref{multCor} to bound $T_m$, where $m = \wh{\mu}$, we need to know the decay behavior of derivatives of $m$. The key estimate that we obtain is
\begin{align*}
\abs{\pd^\alpha_\xi \pd^\beta_\eta\wh{\mu}(\xi,\eta)} 
&\leq C_{\alpha, \beta}( 1 + \min\{\abs{\xi}, \abs{\eta}\}\abs{\sin\theta} )\os{-\frac{d-2}{2}}
\prs{ 1 + \abs{(\xi, \eta)}}\os{-\frac{d-2}{2}},\numberthis\label{eqn7}
\end{align*}
for all multi-indices $\alpha$ and $\beta$. To obtain this estimate we use the following recurrence formula, which can be found in \cite[573]{GrafakosI},
\begin{align*}
	\frac{d}{dt}\prs{t\os{-\nu}J_\nu(t)} = -t\os{-\nu}J_{\nu+1}(t).
\end{align*}
It follows by repeated application of this formula that $\abs{\pd^\alpha_\xi\pd^\beta_\eta \wh{\mu}(\xi,\eta)}$ can be bounded by a finite sum of terms of the form
\begin{align*}
	&~C_{\alpha, \beta}\abs{\xi\os{\gamma}\eta\os{\delta}}(1+\abs{\xi}\abs{\sin\theta})^p(1 + \abs{\eta})^q\int_{0}^1 \abs{ {J_{s}\prs{ 2\pi r \abs{\xi}\abs{\sin\theta} }}}
	\abs{{ J_{t}\prs{2\pi r\frac{\sqrt{3}}{2}\abs{\eta} } }
	}
	\frac{r\dx{r}}{\sqrt{1 - r^2}},
%
\end{align*}
where $s,t > 0$, $\gamma$ and $\delta$ are multi-indices, and
\begin{align*}
	\abs{\xi^\gamma}(\abs{\xi}\abs{\sin\theta})^p = O((1+\abs{\xi}\abs{\sin\theta})\os{-\frac{d-3}{2}}), ~~ \tx{ and }~~ \abs{\eta^\delta}\abs{\eta}^q = O((1+\abs{\eta})\os{-\frac{d-3}{2}}).
\end{align*}
Then one can show using the asymptotic decay of Bessel functions that 
\begin{align*}
	\int_{0}^1 \abs{ {J_{s}\prs{ 2\pi \sqrt{1 - r^2} \abs{\xi}\abs{\sin\theta} }}}
	\abs{{ J_{t}\prs{2\pi \sqrt{1 - r^2}\frac{\sqrt{3}}{2}\abs{\eta} } }
	}\dx{r} &\leq C(1 + \abs{\xi}\abs{\sin\theta})\os{-\frac{1}{2}}(1 + \abs{\eta})\os{-\frac{1}{2}}.
\end{align*}
Applying these estimates and using the fact that $\wh{\mu}(\xi, \eta) = \wh{\mu}(\eta, \xi)$ gives estimate \eqref{eqn7}.

If $\eta = 0$, then we see immediately from \eqref{eqn5} that 
\begin{align*}
	\abs{\wh{\mu}(\xi, 0)} = \abs{\wh{\sigma}_{d-1}(\abs{\xi})} \lesssim (1 + \abs{\xi})\os{-\frac{d-1}{2}}.
\end{align*}
This suggests that the decay behavior of $\wh{\mu}(\xi,\eta)$ is somewhat more complicated than that described by \eqref{eqn7}. However, it also shows that when one of $\abs{\xi}$ or $\abs{\eta}$ is 0 the decay of $\wh{\mu}$ is not much better than that described by \eqref{eqn7}. Similar considerations apply whenever $\eta$ is a positive multiple of $\xi$, as a consequence of Fubini's theorem.


\vskip.25in

\section{$L^p \tm L^p \ra L\os{p/2}$ Bounds Beyond Banach Range}

As described in the introduction, our strategy to obtain upper bounds for operator norms of $T$ as a map from $L^p\tm L^q \ra L^r$ will be as follows. We first decompose the multiplier $m = \wh{\mu}$ dyadically into pieces $m_i$. Then we obtain two bounds for the operator $T_i$ corresponding to the $i$th piece $m_i$ of $m$: $L^2 \ra L^2 \ra L^1$ bounds with norms exponentially decaying in $i$, and $L^1 \tm L^p \ra L\os{\frac{p}{p+1}}$ and $L^p \tm L^1 \ra L\os{\frac{p}{p+1}}$ bounds with norms growing in $i$. Then we interpolate between these bounds as far as possible to get boundedness for the $T_i$ on $L^p \tm L^p \ra L\os{p/2}$ for some range of $p$ such that $p/2 < 1$, with operator norms summable in $i$. We can then bound the operator norm of $T$ on $L^p \tm L^p \ra L\os{p/2}$ by the sum of the operator norms of the $T_i$'s. The key estimates in this procedure are the $L^2 \ra L^2 \ra L^1$ bounds with good operator norms. The following lemma of Grafakos, He, and Slavikova \cite{GrafakosL2} is the essential tool enabling these estimates.

\begin{lemma}\label{lemmaLqMult}
	Let $1 \leq q < 4$ and set $M_q = \floor*{\frac{2d}{4 - q}} + 1$. Let $m$ be a function in $L^q(\eun{2d})\cap C^{M_q}(\eun{2d})$ satisfying
	\begin{equation*}
	\norm{\pd^\alpha m}_{L^\nf} \leq C_0 < \nf
	\end{equation*}
	for all $\abs{\alpha}\leq M_q$. Then there is a constant $D$ depending on $d$ and $q$ such that the bilinear operator $T_m$ with multiplier $m$ satisfies
	\begin{equation*}
	\norm{T_m}_{L^2\tm L^2 \ra L^1} \leq DC_0\os{1 - \frac{q}{4}}\norm{m}_{L^q}\os{\frac{q}{4}}.
	\end{equation*}
\end{lemma}

\noindent We get an immediate corollary to this lemma by estimating $\norm{m}_{L^q}$ for compactly supported $m$.

\begin{corollary}\label{multCor}
	If $m \in L^q(\eun{2d}) \cap C_0^{M_1}(\eun{2d})$, where $M_1 = \floor*{\frac{2d}{3}} + 1$, satisfies
	\begin{equation*}
	\norm{\pd^\alpha m}_{L^\nf} \leq C_0 < \nf,
	\end{equation*}
	for all $\abs{\alpha}\leq M_1$,
	then there is a constant $D$ depending only on $d$ and $q$ such that the bilinear operator $T_m$ with multiplier $m$ satisfies 
	\begin{align*}
		\norm{T_m}_{L^2\tm L^2 \ra L^1} \leq DC_0\abs{\mathrm{sppt}\, m}\os{\frac{1}{4}}.
	\end{align*}
\end{corollary}

\noindent We will decompose the multiplier $m = \wh{\mu}$ into pieces $m_{i,j,k}$, apply this corollary to estimate the operator norm of the operator $T_{i,j,k}$ corresponding to each multiplier $m_{i,j,k}$, then sum the operators $T_{i,j,k}$ to bound the full operator $T$. We define $m_i = \sum_{j,k}m_{i,j,k}$, so that for each fixed $i$ the $m_{i,j,k}$ give a further decomposition of $m_i$. This further decomposition is necessary because $m$ does not have enough uniform decay to apply Corollary \ref{multCor} directly to $m_i$, as discussed in the introduction.


\vskip\subskip

\subsection*{Defining the Decomposition}

We will decompose $m := \wh{\mu}$ using the following three partitions of unity.\\ 

\noindent\textit{1.\ Partition in $\abs{(\xi,\eta)}$}: First we define a partition of unity on a dyadic scale. For this purpose let $\phi_0 \in C_0^\nf(\eun{2d})$ satisfy $\chr{B(0,1)}\leq \phi_0 \leq \chr{B(0,2)}$, and define, for $j \geq 1$, $\phi_j((\xi,\eta)) = \phi_0(2\os{-j}(\xi,\eta)) - \phi_0(2\os{-(j-1)}(\xi,\eta))$. Note that $\phi_j$ is supported in the annulus 
$\{(\xi,\eta) \cl 2\os{(j-1)} \leq \abs{(\xi,\eta)} \leq 2\os{(j+1)}  \}$, and
\begin{equation*}
	\sum_{j=0}^{\nf} \phi_j(\xi,\eta) \equiv 1.
\end{equation*}\slskip

\noindent\textit{2.\ Partition in $\abs{\eta}/\abs{\xi}$}: Now we define a partition of unity regulating the ratio $\abs{\eta}/\abs{\xi}$. First we choose $\psi^* \in C^\nf_0(\rnm)$ with $\chr{[0,1]}\leq \psi^* \leq \chr{[-\epsilon, 1 + \epsilon]}$. Then we define for $j \in \znm$,
\begin{align*}
	\psi^*_j(t) := \frac{\psi^*(t - j)}{\sum_{j\in\znm}\psi^*(t - j)},
\end{align*}
noting that $\psi^*_j$ is supported in $[j-\epsilon, j+1 + \epsilon]$ and $\sum_{j=-\nf}^{\nf}\psi^*_j\equiv 1$.
Finally, we define 
\begin{align*}
	\psi_j(\xi, \eta) := \psi^*_j(\log\abs{\eta} - \log\abs{\xi}) + \psi^*_{-j-1}(\log\abs{\eta} - \log\abs{\xi}),
\end{align*}
so that $\psi_j$ is supported in 
\begin{align*}
	\cbrks{ (\xi,\eta) ~\cl~ 2\os{-\epsilon}2\os{j} \leq \frac{\min\{\abs{\xi},\abs{\eta}\}}{\max\{\abs{\xi},\abs{\eta}\}}\leq 2\os{\epsilon}2\os{j+1} }
\end{align*}
We also define for $j\geq 0$,
\begin{align*}
	\psi^j(\xi,\eta) := \sum_{k=-\nf}^j\psi_k (\xi,\eta),
\end{align*}
which is supported in 
\begin{align*}
\cbrks{ (\xi,\eta) ~\cl \frac{\min\{\abs{\xi},\abs{\eta}\}}{\max\{\abs{\xi},\abs{\eta}\}}\leq 2\os{\epsilon}2\os{j+1} }
\end{align*}
and identically equal to 1 in
\begin{align*}
\cbrks{ (\xi,\eta) ~\cl \frac{\min\{\abs{\xi},\abs{\eta}\}}{\max\{\abs{\xi},\abs{\eta}\}}\leq 2\os{-\epsilon}2\os{j} }.
\end{align*}
In particular this implies that the derivatives of $\psi^j$ vanish except where $\frac{\min\{\abs{\xi},\abs{\eta}\}}{\max\{\abs{\xi},\abs{\eta}\}}\approx 2\os{j}$.\\

\noindent\textit{3.\ Partition in $\abs{\sin\theta}$}: Finally, we define a partition of unity regulating $\abs{\sin\theta}$, where $\theta$ is the angle between $\xi$ and $\eta$. 
First we construct a partition of unity $\{ \rho^*_j \}$ on $\rnm$ on a double-dyadic scale. We let $\rho^*\in \cnm^\nf_0(\rnm)$ with $\chr{[-1,1]}\leq \rho \leq \chr{[-2,2]}$ and we define $\rho^*_j(t) := \rho^*(2\os{-2j}t) - \rho^*(2\os{-2(j-1)}t)$, so that $\rho^*_j$ is supported on the annulus $\{ t ~\cl~ 2\os{2(j-1)} \leq \abs{t} \leq 2\os{2(j+1)} \}$ and $\sum_{j=-\nf}^{\nf}\rho^*_j \equiv 1$, except at 0. 
Now we define, for $j \geq 1$, 
\begin{align*}
	\rho_j(\xi,\eta) = \rho^*_{-j}\prs{ 1 - \prs{ \frac{\xi \cdot \eta}{\abs{\xi}\abs{\eta}} }^2 },
\end{align*}
noting that the argument of $\rho^*_{-j}$ is $(\sin \theta)^2$ and that $\rho_j$ is supported on the set $\{ (\xi,\eta) ~\cl~ 2\os{-j-1}\leq \abs{\sin\theta} \leq 2\os{-j+1} \}$. We also define
\begin{align*}
	\rho_0(\xi,\eta) = \sum_{j=0}^{\nf}\rho^*_{j}\prs{ 1 - \prs{ \frac{\xi \cdot \eta}{\abs{\xi}\abs{\eta}} }^2 },
\end{align*}
so that $\rho_0(\xi,\eta)$ is supported on the set $\{ (\xi,\eta) ~\cl \abs{\sin\theta} \geq 1/2 \}$. We define for $j \geq 0$
\begin{align*}
	\rho^j(\xi, \eta) := \sum_{k=\nf}^{-j}\rho^*_{k}\prs{ 1 - \prs{ \frac{\xi \cdot \eta}{\abs{\xi}\abs{\eta}} }^2 },
\end{align*}
so that $\rho^j$ is supported on $\{ (\xi, \eta) ~\cl \abs{\sin\theta} \leq 2\os{-j+1} \}$ and is identically equal to 1 on $\{ (\xi, \eta) ~\cl \abs{\sin\theta} \leq 2\os{-j} \}$. Note that this implies that the derivatives of $\rho^j$ vanish on $\{ (\xi, \eta) ~\cl \abs{\sin\theta} \leq 2\os{-j} \}$.\\

We now apply these three partitions to $m$. First we fix $i \in \{0, 1, 2, \ldots\}$. Then for $0\leq j < i$ and $0 \leq k < {\floor{\frac{i-j}{2}}  }$ we define
\begin{align*}
	m_{i,j,k}:= m\phi_i\psi_j\rho_k
\end{align*}
and for $0 \leq j \leq i$ we define 
$$m_{i,j,\floor{\frac{i-j}{2}}}:= m\phi_i\psi_j\rho\os{\floor{\frac{i-j}{2}}},$$ noting that $\rho^0 \equiv 1$. 
Recall that on the support of $m_{i,j,k}$, $\xi$ and $\eta$ satisfy
\begin{align*}
2\os{2(i-1)} \leq \abs{\xi}^2 + \abs{\eta}^2 \leq 2\os{2(i+1)}
~~ \tx{ and } ~~~ 2\os{-\epsilon}2\os{j} \leq \frac{\min\{\abs{\xi},\abs{\eta}\}}{\max\{\abs{\xi},\abs{\eta}\}}\leq 2\os{\epsilon}2\os{j+1},
\end{align*}
and hence
\begin{align*}
2\os{i-j +2} \geq \min\{ \abs{\xi},\abs{\eta} \} \geq 2\os{i-j-3}.
\end{align*}

We define the multiplier
\begin{align*}
m_i:= \sum_{j=0}^{i}\sum_{k=0}^{\floor{\frac{i-j}{2}}}m_{i,j,k},
\end{align*}
and the corresponding operators
\begin{align*}
T_i(f,g)(x) := \ft\finv\prs{ m_i\wh{f}\otimes\wh{g} }(x,x) ~~ \tx{ and } ~~ T_{i,j,k}(f,g)(x) := \ft\finv\prs{ m_{i,j,k}\wh{f}\otimes\wh{g} }(x,x),
\end{align*}
for Schwartz functions $f$ and $g$. We will estimate right-hand side of
\begin{align*}
\norm{T_i}_{L^2\tm L^2 \ra L^1}\leq \sum_{j=0}^{i}\sum_{k=0}^{\floor{\frac{i-j}{2}}}\norm{T_{i,j,k}}_{L^2\tm L^2 \ra L^1},
\end{align*}
in order to obtain bounds for $\norm{T_i}_{L^2\tm L^2 \ra L^1}$ which are exponentially decaying in $i$.


\newpage


\subsection*{Computing derivatives of $m_{i,j,k}$}
In order to apply Corollary \ref{multCor}, we need $L^\nf$ bounds for derivatives of the multipliers $m_{i,j,k}$. 
We have seen above that all derivatives of $m$ satisfy the same decay estimates. Hence when we take derivatives of $m_{i,j,k}$ we have to consider the size of the derivatives of the cutoffs $\phi_i$, $\psi_j$, and $\rho_k$ on the support of $m_{i,j,k}$. We claim that these derivatives are uniformly bounded. In estimating these derivatives we use the following: 

\begin{heuristic}
If $x\os{\alpha}\abs{x}^p$ is $\lesssim \abs{x}^k$, then
 $\pd_{x_i}\prs{x\os{\alpha}\abs{x}^p}$ is $\lesssim \abs{x}\os{k-1}$, unless it is 0. 
\end{heuristic}
\noindent Using this heuristic and induction, one can estimate that on the support of $m_{i,j,k}$, 
\begin{align*}
	\norm{\pd^\alpha_{\xi}\pd^\beta_{\eta}\phi_i}_\nf \leq C,
\end{align*}
that
\begin{align*}
	\norm{\pd^\alpha_{\xi}\pd^\beta_{\eta}\psi_j}_\nf \lesssim \abs{\xi}\os{-\abs{\alpha}}\abs{\eta}\os{-\abs{\beta}} \lesssim 2\os{-(\abs{\alpha} + \abs{\beta})(i-j)},
\end{align*}
and
\begin{align*}
	\norm{\pd^\alpha_{\xi}\pd^\beta_{\eta}\rho_k}_\nf \lesssim 2\os{2k(\abs{\alpha}+\abs{\beta})}\abs{\xi}\os{-\abs{\alpha}}\abs{\eta}\os{-\abs{\beta}} 
	 \lesssim 
	 2\os{-(\abs{\alpha} + \abs{\beta})(i-j - 2k)}.
\end{align*}
We see that by the choice of the indices $i,j,k$ included in our decomposition these derivatives are uniformly bounded on the support of $m_{i,j,k}$. 


\vskip\subskip

\subsection*{Volume of $\spt m_{i,j,k}$}

The support of $m_{i,j,k}$ is contained in
\begin{align*}
	S_{i,j,k} = \cbrks{ (\xi,\eta) ~\cl~ \abs{\xi}\leq 2\cdot 2\os{i},~ \abs{\eta} \leq 4\cdot 2\os{i-j}, ~ 2\os{-k-1}\leq \abs{\sin\theta} \leq 2\os{-k+1} }.
\end{align*}
We can estimate the volume of this set using radial coordinates and the slicing formula \eqref{eqn3}. Defining $A_{k} := \{ (\xi,\eta) ~\cl~ 2\os{-k-1}\leq \abs{\sin\theta} \leq 2\os{-k+1} \}$, where $\theta$ is the angle between $\xi$ and $\eta$, 
\begin{align*}
	\abs{S_{i,j,k}}&= \int_{0}^{2\cdot2\os{i}}\int_{0}^{2\os{i-j+2}}\int_{\sph{d-1}}\int_{\sph{d-1}}\chr{ A_k}(r\omega, s\nu)\dx{\sigma\prs{\omega}}\dx{\sigma(\nu)}r\os{d-1}\dx{r}\,s\os{d-1}\dx{s}\nlspace
	&= C_d2\os{id}2\os{(i-j)d}\int_{\sph{d-1}}\int_{\sph{d-1}}\chr{A_k}(\omega, \nu)\dx{\sigma\prs{\omega}}\dx{\sigma(\nu)}\nlspace
	&= C_d2\os{id}2\os{(i-j)d}\int_{\sph{d-1}}\int_{\sph{d-1}}\chr{A_k}(\omega, e_1)\dx{\sigma\prs{\omega}}\dx{\sigma(\nu)}\nlspace
	&= C_d2\os{id}2\os{(i-j)d}\int_{\sph{d-1}}\chr{}\cbrks{ 2\os{-k - 1}\leq \sqrt{1 - \omega_1^2} \leq 2\os{-k+1} }\dx{\sigma}\nlspace
	&= C_d2\os{id}2\os{(i-j)d} \int_{0}^{1} \chr{}\cbrks{ 2\os{-k - 1}\leq r \leq 2\os{-k+1} } \frac{r\os{d-2}}{\sqrt{1 - r^2}}\dx{r}\nlspace
	& \approx_d 2\os{id}2\os{(i-j)d}2\os{-k(d-1)}.
\end{align*}
One can show that $\abs{\spt m_{i,j,k}}\approx \abs{S_{i,j,k}}$ with dimensional constants. 

%
%
%
%


\vskip\subskip

\subsection*{$L^2 \tm L^2 \ra L^1$ bound for $T_i$}

Using Corollary \ref{multCor} and \eqref{eqn7} and the fact that the derivatives of the cutoff functions are uniformly bounded in $i,j$, and $k$, for $0\leq j \leq i$ and $0 \leq k < \max\cbrks{\floor{\frac{i-j}{2}},1}$ we get the bound 
\begin{align*}
	\norm{T_{i,j,k}}_{L^2\tm L^2 \ra L^1} \leq C 2\os{-i\frac{d-2}{2}}\cdot\sbrks{2\os{i-j-k}}\os{-\frac{d-2}{2}}\cdot 2\os{(2i - j - k)\frac{d}{4}}2\os{\frac{k}{4}}.
\end{align*}
For $0 \leq j < i - 1$ and $k = \floor{\frac{i-j}{2}}$ the following bound is derived in Appendix \ref{a1}:
\begin{align*}
	\norm{T_{i,j,k}}_{L^2\tm L^2 \ra L^1} &\leq C2\os{i\frac{-3d + 20}{16}} \cdot 2\os{j\frac{-d - 4}{16}}.
\end{align*}
Thus summing over $j$ and $k$ gives the bound
\begin{align*}
	\norm{T_i}_{L^2\tm L^2 \ra L^1} &\leq C2\os{-i(\frac{d}{2}-2)}\cdot 
		\sum_{j=0}^{i} 2\os{j\prs{ \frac{d}{4} - 1 }}
		\sum_{k=0}^{\max\cbrks{\floor{\frac{i-j}{2}}-1, 0}} 2\os{k\prs{ \frac{d}{4} - \frac{3}{4} }}
		+ C\sum_{j=0}^{i-2}2\os{i\frac{-3d + 20}{16}} \cdot 2\os{j\frac{-d - 4}{16}}
	\nlspace
	&\leq C2\os{-i\prs{\frac{d}{2}-2}}\cdot 
		\sum_{j=0}^{i}\, 2\os{j\prs{ \frac{d}{4} - 1 }}\cdot 2\os{{\prs{\frac{i-j}{2}}}\prs{ \frac{d}{4} - \frac{3}{4} }}
		+ C2\os{i\frac{-3d + 20}{16}}
	\nlspace
	&\leq C 2\os{-i\prs{ \frac{d}{4} - 1 }} + C2\os{i\frac{-3d + 20}{16}}\nlspace
	&\leq C2\os{i\frac{-3d + 20}{16}}.
\end{align*}
which is exponentially decreasing in $i$ for $d \geq 7$. 


\vskip\subskip

\subsection*{$L^1 \tm L^p \ra L\os{\frac{p}{p+1}}$ bound for $T_i$}

We use the following inequality, which can be found in \cite[480]{GrafakosI}. For $N > M > d$,
\begin{align*}
	\int_{\sph{d-1}}\frac{2\os{dj}}{\prs{1 + 2\os{j}\abs{x - y}}^N} \leq \frac{C_{M,N}2\os{j}}{\prs{1 + \abs{x}}^M}.
\end{align*}
Now observe that for $x = (x_1, x_2)\in \eun{2d}$, using the fact that $\phi$ is a Schwartz function, we have
\begin{align*}
	\abs{(\widecheck{\phi}_i\star \mu)(x)} &\leq C_N\int_{SO(d)}\frac{2\os{2di}}{\prs{ 1 + 2^i\abs{(x_1 - Re_1, x_2 - R((1/2)e_1 + (\sqrt{3}/2)e_2)} }^{2N}}\dx{R}\nlspace
	&\leq 2\os{di}\prs{ \int_{\sph{d-1}} \frac{2\os{di}}{(1 + 2\os{i}\abs{ x_1 - \omega })^{2N}}\dx{\sigma\omega} }\os{\frac{1}{2}}\prs{ \int_{\sph{d-1}} \frac{2\os{di}}{(1 + 2\os{i}\abs{ x_2 - \omega })^{2N}}\dx{\sigma\omega} }\os{\frac{1}{2}}\nlspace
	&\leq \frac{C_{M,N}2\os{(d+1)i}}{ (1 + \abs{x_1})^M(1 + \abs{x_2})^M }.
\end{align*}
Hence we have that
\begin{align*}
	\norm{T_i(f,g)}_{L\os{\frac{p}{p+1}}} &\leq C2\os{(d+1)i}\norm{ (1 + \abs{x})\os{-M}\star \abs{f} }_{L^{1}}\norm{ (1 + \abs{x})\os{-M}\star \abs{g} }_{L^{p}}\nlspace
	&\leq C2\os{(d+1)i}\norm{ f }_{L^{1}}\norm{g }_{L^{p}},
\end{align*}
by Young's inequality. Analogously we get
\begin{align*}
	\norm{T_i(f,g)}_{L\os{\frac{p}{p+1}}} &\leq C2\os{(d+1)i}\norm{ f }_{L^{p}}\norm{g }_{L^{1}}.
\end{align*}



\vskip.25in

\section{Interpolation}

\vskip.125in

We first obtain $L^p \tm L^q \ra L^r$ bounds in the range of indices $(p,q,r)$ for which $1/p + 1/q \leq 1$. We refer to this range of indices as the ``Banach range'' because $1/r \leq 1/p + 1/q \leq 1$ must hold \cite[Prop.\ 7.1.5]{GrafakosII}. 
We can majorize $T$ as follows:

\begin{align*}
\abs{T(f,g)(x)} &= \abs{\int_{SO(d)}f\prs{x - Re_1^d}\int_{SO(d-1)}g\prs{x- \frac{1}{2}Re_1^d - { \frac{\sqrt{3}}{2}R(0,R'e_1^{d-1})}}\dx{R'}\dx{R}}\nlspace
%
%
&\leq C\norm{g}_{L^\nf}S_1(\abs{f})(x),
\end{align*}
where $S_1$ is the linear spherical averaging operator with radius 1 acting on functions on $\eun{d}$.
We know that $S_1$ is bounded from $L^p$ into $L^q$ whenever $\prs{ 1/p, 1/q }$ lies in convex hull $\conv V$ of the set of points
\begin{align*}
V:= \cbrks{(0,0), (1,1), \prs{\frac{d}{d+1}, \frac{1}{d+1}}}.
\end{align*}
Hence $T$ is bounded from 
\begin{align*}
L^p(\eun{d}) \tm L^\nf(\eun{d}) \ra L^q(\eun{d})
\end{align*}
whenever $(1/p, 1/q)$ lies in $\conv V$. 
By the symmetry of the manifold, we can exchange the roles of $f$ and $g$ in the above integral to also get boundedness from 
\begin{align*}
L^\nf(\eun{d}) \tm L^p(\eun{d}) \ra L^q(\eun{d})
\end{align*}
whenever $(1/p, 1/q)$ lies in $\conv V$. 
We can now apply bilinear interpolation to see that $T$ is bounded from
\begin{align*}
L^p(\eun{d}) \tm L^q(\eun{d}) \ra L^r(\eun{d})
\end{align*}
whenever $(1/p, 1/q)$ lies in the convex hull of the set 
\begin{align*}
V_2:= \cbrks{ (0,0), (1,0), (0,1) },
\end{align*}
and we can push $1/r$ some amount below $1/p + 1/q$, depending on the choice of $p, q$, using the fact that $S_1$ is $L^p$-improving.
One can also obtain these Banach range bounds via Minkowski's integral inequality.


Interpolating between the bounds obtained in the last section above gives that
\begin{align*}
	\norm{T_i(f,g)}_{L\os{p/2}} &\leq C2\os{(d+1)i}\norm{f}_{L^p}\norm{g}_{L^p},
\end{align*}
for $ 1 < p < \nf$. We now want to interpolate against the $L^2 \tm L^2 \ra L^1$ bound as far as possible while still maintaining the summability in $i$ of the norms of the operators ${T_i}_{}$. Doing so, we find that $T$  maps $L^p\tm L^p \ra L\os{p/2}$
with operator norms summable in $i$ for $$p > p_d := \frac{19d - 4}{11d - 12}.$$ Finally, we can interpolate these bounds against those for exponents in the Banach range to get boundedness on the interior of the convex hull of the set of indices $\{(0,0), (1,0), (0,1), (1/p_d, 1/p_d)\}$, for $d \geq 7$, as displayed in Figure \ref{fig1}.


\vskip.25in

\section{Maximal Operator}
\newcommand{\tmax}{\mc{T}}

\vskip.125in

We define the maximal operator that arises from allowing the radius of triangles to vary in the definition of $T$:
\begin{align*}
	\mc{T}(f,g)(x) := \sup_{t > 0}\abs{ \int_M f(x - tu)g(x - tv) \dx{\mu(u,v)}}.
\end{align*}
We get a range of bounds for this operator immediately by the majorizations
\begin{align*}
	T(f,g)(x) \leq \norm{g}_{L^\nf}\sup_{t > 0}\abs{ S_t(f)(x) }, ~~\text{ and }~~~~~ T(f,g)(x) \leq \norm{f}_{L^\nf}\sup_{t > 0}\abs{ S_t(g)(x) }
\end{align*}
where
\begin{align*}
	S_t(h)(x) = \int_{\sph{d-1}}\abs{h(x-t\omega)}\dx{\sigma_{d-1}(\omega)}
\end{align*}
is the spherical averaging operator with radius $t > 0$. By the well-known results of Stein and Bourgain, $\mc M:= \sup_{t > 0}\abs{S_t}$ is bounded from $L^p \ra L^p$ for $p \in (\frac{d}{d-1},\nf]$, for $d\geq 2$. We can now interpolate between these bounds to show that $\tmax$ is bounded for $(\frac{1}{p},\frac{1}{q})$ in the closed convex hull of the set of points $\{ (0,0), (0,\frac{d-1}{d}), (\frac{d-1}{d}, 0) \}$, minus the upper-right boundary. By contrast, in \cite{Barrio} the authors show that the procedure of taking the $L^\nf$ norm of one function and applying a linear multiplier theorem gives boundedness on the full Banach range for the maximal bilinear spherical averaging operator. 

We can use an example similar to the one used in \cite{Barrio} to establish lower bounds on the indices $p,q$ such that $\mc T$ is bounded on $L^p(\eun{d}) \tm L^q(\eun{d})$.
We set 
\begin{equation*}
f(x) := \frac{1}{\abs{x}\os{\frac{d}{p}}(-\log\abs{x})\os{\frac{2}{p}}}\chr{\{ \abs{x} \leq 1/8 \}}(x) ~~\text{ and }~~
g(x) := \frac{1}{\abs{x}\os{\frac{d}{q}}(\log\abs{x})\os{\frac{2}{q}}}\chr{\{ \abs{x} \geq 8 \}}(x).
\end{equation*}
Then $f \in L^p(\eun{d})$ and $g \in L^q(\eun{d})$, and by taking $t = R$ we can estimate that
\begin{align*}
	\mc T(f,g)(x) &\gtrsim 
	\begin{cases}
	\frac{ \chr{\{ \abs{x} \geq 7/8 \}}(x) }{\abs{x}\os{ d\prs{1 + \frac{1}{q} - \frac{1}{p}} + \epsilon} \prs{\log \abs{x}}\os{\frac{2}{q}} } & \text{ if $p \geq \frac{d}{d-1}$}\nlspace
	\nf & \text{ if $p < \frac{d}{d-1}$}.
	\end{cases}
\end{align*}
Now by using the assumption that $\mc T(f,g)(x) \in L^r(\eun{d})$ and interchanging $f$ and $g$ we can deduce that both $p$ and $q$ are $> \frac{d}{d-1}$.
It follows that the maximal region of boundedness for $\mc T$ on $L^p(\eun{d}) \tm L^q(\eun{d})$ in terms of the exponents $(1/p, 1/q)$ is contained inside the square shown in Figure \ref{fig2}.\\[0mm]
\begin{center}
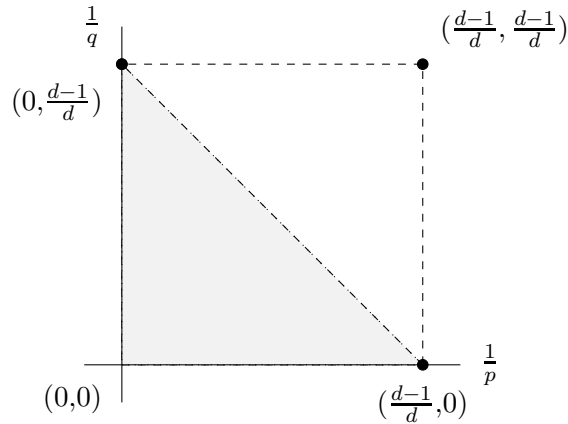

	\begin{tikzpicture}[scale=0.5]
	\draw[style = dashed] (0,0) -- (8,0) -- (8,8) -- (0,8) -- cycle;
	\draw[fill=gray!10, style = dotted] (0,0) -- (8,0) -- (0,8) -- cycle;
	\draw[color=black] (-1,0) -- (9,0);
	\draw[color=black] (0,-1) -- (0,9);
	\draw[color=black, style=dashed] (8,0) -- (0,8);
	
	
	\node [fill, draw, circle, minimum width=4pt, inner sep=0pt] at (8,8) {};
	\node[above right=1pt of {(8,8)}, outer sep=2pt] {\small$(\frac{d-1}{d},\frac{d-1}{d})$};
	
	\node [fill, draw, circle, minimum width=4pt, inner sep=0pt] at (8,0) {};
	\node[below=2pt of {(8,0)}, outer sep=2pt] {\small($\frac{d-1}{d}$,0)};
	
	\node [fill, draw, circle, minimum width=4pt, inner sep=0pt] at (0,8) {};
	\node[below left=2pt of {(0,8)}, outer sep=2pt] {\small(0,$\frac{d-1}{d}$)};
	
	\node [fill, draw, circle, minimum width=4pt, inner sep=0pt] at (0,8) {};
	\node[below left=2pt of {(0,0)}, outer sep=2pt] {\small(0,0)};
	
	\node[right=2pt of {(9,0)}, outer sep=2pt] {$\frac{1}{p}$};
	
	\node[left=2pt of {(0,9)}, outer sep=2pt] {$\frac{1}{q}$};
	\end{tikzpicture}
	\captionof{figure}{\color{black} Maximum region of boundedness for $\mc T$ and region obtained by interpolation. }\label{fig2}
\end{center}

\vspace*{5mm}

\noindent It remains to determine whether $\mc{T}$ is bounded on $L^p(\eun{d}) \tm L^q(\eun{d})$ for any pair $(1/p, 1/q)$ in the upper-right half of this square. 

\appendix

\vskip.25in

\section{Derivation of bound for $\norms{T_{i,j,\floor{\frac{i-j}{2}}}}_{L^2\tm L^2 \ra L^1}$}\label{a1}

\vskip.125in

\newcommand{\floors}[1]{\lfloor#1\rfloor}

We need to estimate the operator norm of the multiplier $m_{i,j,\floors{\frac{i - j}{2}}}(\xi,\eta)$. This multiplier is supported where $|\xi| \sim 2\os{i}$, $|\eta|\sim 2\os{i-j}$, and $0 \leq \abs{\sin \theta} \leq 2\os{-\floor{\frac{i - j}{2}}}$. The best estimate we have for the $L^\nf$ norm of $m_{i,j,\floors{\frac{i - j}{2}}}$ and its derivatives that holds on all of $\spt(m_{i,j,\floors{\frac{i - j}{2}}})$ is 
\[
\abs{\pd^\alpha m_{i,j,\floors{\frac{i - j}{2}}}(\xi, \eta)} \leq C_\alpha 2\os{-i\frac{d-2}{2}}.
\]
Hence by Lemma 3,
\[
\norms{T_{m_{i,j,\floors{\frac{i - j}{2}}}}}_{L^2\tm L^2 \ra L^1} \leq C2\os{-i\frac{d-2}{2}\cdot\frac{3}{4}}\norms{m_{i,j,\floors{\frac{i - j}{2}}}}_{L^1}\os{\frac{1}{4}}.\numberthis\label{e1}
\]
We can now estimate
\begin{align*}
\norms{m_{i,j,\floors{\frac{i - j}{2}}}}_{L^1} &\leq C2\os{-i\frac{d-2}{2}}\int_{\spt(m_{i,j,\floors{\frac{i - j}{2}}})}\prs{1 + 2\os{i-j}|\sin\theta|}\os{-\frac{d-2}{2}}  \dx{(\xi,\eta)},
\end{align*}
and we can bound the integral by a constant times
\begin{align*}
\int_{\substack{|\xi|\sim 2^i, |\eta|\sim 2\os{i-j},\\ 0\leq |\sin\theta| \leq 2\os{-(i-j)}}}\dx{(\xi,\eta)} 
+ 2\os{-(i-j)\frac{d-2}{2}}\int_{\substack{|\xi|\sim 2^i, |\eta|\sim 2\os{i-j},\\ 
		2\os{-(i-j)} \leq |\sin\theta| \leq 2\os{-(i-j)/2} }}|\sin\theta|\os{-\frac{d-2}{2}}\dx{(\xi,\eta)}.
\end{align*}
As in the section \textbf{Volume of $\spt m_{i,j,k}$} above, the first integral is 
\[
\lesssim 2\os{id}2\os{(i-j)d}2\os{-(i-j)(d-1)}.
\]
The contribution to the operator norm from the first integral will be negligible compared to the contribution from the second integral. One can now proceed similarly to estimate the second integral by
\begin{align*}
&\int_{\frac{1}{C}2\os{i}}^{C2\os{i}}\int_{\frac{1}{C}2\os{i-j}}^{C2\os{i-j}}
\int_{\sph{d-1}}\int_{\sph{d-1}}
\chr{}\{ 2\os{-(i-j)} \leq |\sin\theta(r\omega, s\nu)| \leq 2\os{-(i-j)/2} \}\nlspace
&\hspace*{5cm}\cdot|\sin\theta(r\omega, s\nu)|\os{-\frac{d-2}{2}}
\dx{\sigma(\omega)}\dx{\sigma(\nu)}
r\os{d-1}\dx{r}s\os{d-1}\dx{s}\nlspace
%
%
&\quad\leq C_d 2\os{id}2\os{(i-j)d}\int_{2\os{-(i-j)}}^{2\os{-(i-j)/2}}t\os{-\frac{d-2}{2}}\frac{t\os{d-2}}{\sqrt{1 - t^2}}\dx{t}
~\approx_d~ 2\os{id}2\os{(i-j)d}2\os{-(i-j)\frac{d}{4}}.
\end{align*}
Now inserting these estimates back into \eqref{e1} and collecting terms, we have that 
\begin{align*}
\norms{T_{m_{i,j,\floors{\frac{i - j}{2}}}}}_{L^2\tm L^2 \ra L^1} &\leq C2\os{-i\frac{d-2}{2}}\cdot 2\os{(2i-j)\frac{d}{4}}\cdot 2\os{-(i - j)\frac{3d - 4}{16}} \nlspace
&= C2\os{i\frac{-3d + 20}{16}} \cdot 2\os{j\frac{-d - 4}{16}}. 
\end{align*}

%
%
%
%
%
%
%
%

\end{document}